\documentclass[14pt]{amsart}
\usepackage{latexsym,enumerate}
\usepackage{amsmath,amsthm,amsopn,amstext,amscd,amsfonts,amssymb}
\usepackage[latin1]{inputenc}
\usepackage{graphicx}
\usepackage[active]{srcltx}

\newcommand{\dem}{\noindent\bf Proof.\rm $\:$ }

\newcommand{\limit}{\lim\limits}

\renewcommand{\b }{\beta }

\newcommand{\e }{\varepsilon }

\newcommand{\cqd}{{\unskip\nobreak\hfil\penalty50
        \hskip2em\hbox{}\nobreak\hfil\mbox{\rule{1ex}{1ex} \qquad}
        \parfillskip=0pt \finalhyphendemerits=0\par\medskip}}

\newtheorem{Theorem}{Theorem}[section]
\newtheorem{Corollary}[Theorem]{Corollary}

\newtheorem{Lemma}[Theorem]{Lemma}
\newtheorem{Proposition}[Theorem]{Proposition}

\newtheorem{remark}[Theorem]{Remark}
\date{}
\catcode`@=11 \@addtoreset{equation}{section}
\renewcommand\theequation{\thesection.\@arabic\c@equation}
\catcode`@=12

\begin{document}

\title{Global analysis of an infection age SEI model with a large class of nonlinear incidence rates.}

\author{ Sofiane Bentout}
\author{Tarik Mohamed Touaoula}

\address{Department of Mathematics, University
Aboubekr Belka\"{\i}d,  Tlemcen 13000, Algeria.}
\email{{\tt{bentoutsofiane@gmail.com}}}
\email{{\tt{tarik.touaoula@mail.univ-tlemcen.dz}}}

\begin{abstract}
We propose and investigate an SEI infection's age model with a general class of nonlinear incidence rates. We give a necessary and sufficient condition for global asymptotic stability of the free-equilibrium related to the basic reproduction number. By using  Lyapunov functionals, we show the global asymptotic stability of the  endemic equilibrium whenever it exists.
\end{abstract}


\maketitle




\section{Introduction}\label{sec:1}



In the past decades, many epidemiological models with infection age have been studied. At our knowledge, The paper \cite{Thieme1} is a leading work that showed the impact of the infection's age  in such models.\\

Another important aspect in disease modeling, is the incidence rate, that is, a newly infected individuals per unit of time.

A bilinear form of incidence rate (action mass type) such as, $SI,$  or $S\int_0^{\infty}\beta(a)i(.,a)da,$ ($\beta(a)$ being the transmission rate with respect to the age of infection and $i$ being the density of infection individuals) is frequently used; this is often considered to characterize the fact that the contact number between susceptible and infective is proportional to the product of both sub-populations, see for instance (\cite{Brower}, \cite{Huang2}, \cite{Korobeinikov}, \cite{Song}, \cite{Ma}, \cite{McCluskey}, \cite{McCluskey1}, \cite{McCluskey2}, \cite{Rost},   \cite{Takeuchi}) and the references therein.

However many authors considered that the bilinear functional is not always realistic; and one have to consider some other forms than the bilinear one to ensure a good qualitative description of the disease dynamics.

Some of the works considering a more general incidence rate functional, are those of  Feng and Thieme \cite{feng}, \cite{Feng} and Korobeinikov et al.  \cite{Korobeinikov1},\cite{Korobeinikov2}, where the incidence was modeled in the nonlinear form $f(S,I),$ and some global properties of $SIR$ and $SIRS$ epidemic models was established.

Delay epidemic models with nonlinear incidence were studied in \cite{Xu}, \cite{Zhao}, nevertheless, only global stability of the disease free equilibrium was proved. Later Huang et al. \cite{Huang}, \cite{Huang1}, considered the more general nonlinear incidence rate of the form $f(S(t),I(t-\tau));$ In these both papers, global stability of endemic and disease free equilibria was established by constructing suitable Lyapunov functionals.

 Concerning the infection's age SIR model with mass action type, Magal et al. \cite{Magal} employed a convenient Lyapunov functional in order to obtain the asymptotic stability of the endemic equilibrium. Very recently in \cite{bentout}, the infection's SIR model with general nonlinear incidence is investigated.

\noindent  In the context of infection's age SEI model, Rost et al.  \cite{Rost} and McCluskey \cite{McCluskey},  investigated the following model
\begin{equation*}
\left \{
\begin{array}{lll}
 S'(t)=\Lambda-\beta S(t)\int_0^{\infty}k(a)i(t,a)da-dS(t)\\\\
  E'(t)=\beta S(t)\int_0^{\infty}k(a)i(t,a)da-(\mu+d)E(t) \\\\
  i_t(t,a)+i_a(t,a)=-(d+\delta+r)i(t,a),\\\\
R'(t)=r I(t)-d R(t).
\end{array}
\right.
\end{equation*}
subject to the following boundary condition
$$i(t,0)=\mu E(t).$$

 The authors in \cite{Rost} proposed and analyzed this model. More precisely, they proved the global stability of the free equilibrium, and the local stability of the endemic equilibrium. In addition they showed that the disease is always present whenever the basic reproduction number $R_0$ is bigger than one. Right after Clusckey  in \cite{McCluskey}, proved, by using a Lyapunov functional, that the endemic equilibrium is globally asymptotically stable whenever it exists.

Motivated by  these previous works, we propose in this paper an infection's age SEI  model with  a general nonlinear incidence rate :

\begin{equation} \label{A1}
\left \{
\begin{array}{lll}
S'(t)=A-\mu S(t)-f\big(S(t),J(t)\big), \;\;\ t \geq 0,\\\\
E'(t)=f\big(S(t),J(t)\big)-(\mu+\alpha)E(t),\\\\
 i_t(t,a)+i_a(t,a)=-(\mu+\gamma(a)) i(t,a), \\\\
 J(t)=\int_{0}^{\infty}\beta(a)i(t,a)da,
\end{array}
\right.
\end{equation}
with the boundary and initial conditions
\begin{equation*}
\left \{
\begin{array}{lll}
i(t,0)=\alpha E(t),\;\ t>0,\\\\
S(0)=S_0\geq 0, \\\\
i(0,.)=i_0(.),
\end{array}
\right.
\end{equation*}

where $\b(a),$  represents the transmission coefficient and $\gamma(a)$ represents the recovery rate of infected individuals, both considered at age $a.$ The parameter $A$ is the entering flux into the susceptible class $(S)$ and $\mu$ is the mortality rate of the population.

\noindent We will consider in the whole paper that the function $\beta $ is an integrable positive function. The function $\gamma$ belong to the set $L_{+}^{\infty}(\mathbb{R}^{+})$, the  nonnegative cone of $L^{\infty}(\mathbb{R}^{+})$. The parameters $A$ and $\mu$ are supposed to be positive.

 We suppose that the incidence function $f$ satisfies :

  $f(S,.)$ is strictly increasing for $S> 0$ and   $f(.,J)$ is strictly increasing for $J>0.$ Moreover  $f(0,J)=f(S,0)=0$  for all $S, J\geq 0.$  \\\\
\noindent The function  $\dfrac{\partial f}{\partial J}(.,0)$ is continuous positive on every compact set $K$.  \\\\
\noindent The function $f$ is locally Lipschitz continuous in $S$ and $J,$ with a Lipschitz constant $L>0,$ i.e. for every $C>0$ there exists some $L:=L_C>0$ such that
  \begin{equation}\label{Lip}
  |f(S_2,J_2)-f(S_1,J_1)|\leq L(|S_2-S_1|+|J_2-J_1|),
  \end{equation}
  whenever $0\leq S_2,S_1, J_2, J_1\leq C.$

\noindent The disease splits the population ($N$) into susceptible, infective, exposed and removed individuals, that is  $N(t)=S(t)+E(t)+I(t)+R(t)$ with $I(t)=\int_0^{\infty}i(t,a)da$.  The removed class is modeled as
\begin{equation}\label{R}
R'(t)=\int_0^{\infty}\gamma(a)i(t,a)da-\mu R(t).
\end{equation}

The total population $N(t)$ satisfies the following ordinary differential equation :
\begin{equation*}
N'(t)=A-\mu N(t).
\end{equation*}
It is obvious that $N(t)$ converges to $\dfrac{A}{\mu},$ when $t\rightarrow \infty$ hence, the equation of $R$ in (\ref{R}) can be omitted.\\

Throughout this paper, we denote $$\bar{N}=\dfrac{A}{\mu},$$ and

$$\digamma(a)=e^{-\int_0^{a}\gamma(\sigma)d\sigma}.$$

\noindent Notice that by integrating the equation of $i$ in (\ref{A1}) along the characteristic lines $t-a=constant,$ the solution $(S,i)$ satisfies,
 \begin{equation}\label{phi}
i(t,a)=\left\{
\begin{array}{lll}
\alpha E(t-a)e^{-\mu a}\digamma(a), \;\;\ t>a\geq 0,\\\\
i_0(a-t)e^{-\mu t} \dfrac{\digamma(a)}{\digamma(a-t)}, \;\;\ a>t\geq 0.
\end{array}
\right.
\end{equation}
and
\begin{equation}\label{phi1}
\left\{
\begin{array}{lll}
S'(t)=A-\mu S(t)-f(S(t),J(t))\\\\
E'(t)=f(S(t),J(t))-(\mu+\alpha)E(t),\\\\
S(0)=S_0,\;\ E(0)=E_0.
\end{array}
\right.
\end{equation}
with
\begin{equation*}
 J(t)=\int_0^{\infty}\b(a)i(t,a)da.
\end{equation*}

\noindent  Hence, the problem (\ref{A1}) may be rewritten as
\begin{equation} \label{A2}
\left \{
\begin{array}{lll}
S'(t)=A-\mu S(t)-f\big(S(t),J(t)\big), \,\\\\
E'(t)=f\big(S(t),J(t)\big)-(\mu+\alpha)E(t),\\\\
J(t)=\alpha\int_0^{\infty}\beta(a)e^{-\mu a}\digamma(a)E(t-a)da,
\end{array}
\right.
\end{equation}
with the initial conditions
\begin{equation} \label{A22}
\left \{
\begin{array}{lll}
S(0)=S_0,\\\\
E(t)=\phi(t), \;\ \mbox{for}\;\ t\leq 0.
\end{array}
\right.
\end{equation}

\noindent We denote
\begin{equation*}
C_{\Delta}:=\{\phi:(-\infty, 0]\rightarrow \mathbb{R}, \;\ \phi(s) e^{\Delta s}\;\ \mbox{is bounded and uniformly continuous on}\;\ (-\infty, 0]\}.
\end{equation*}
with $0<\Delta<\mu+\inf(\gamma).$\\
$C_{\Delta}$ is a Banach space endowed with the norm
\begin{equation*}
||\phi||=\sup_{s\leq 0}|\phi(s)e^{\Delta s}|.
\end{equation*}

Let $E_t$ denote the state of the solution $E(t)$ at time $t$; i.e. $E_t(\theta)=E(t+\theta),$ where $\theta\leq 0.$ We are interested only in the nonnegative solutions, the corresponding cone of non-negative functions in $C_{\Delta}$ is denoted by $Y;$ i.e.
\begin{equation*}
Y:=\{\phi \in C_{\Delta} : \phi(\theta)\geq 0\;\ \mbox{for}\;\ \theta\leq 0\}.
\end{equation*}

We assume that $(S_0,\phi) \in \mathbb{R}^{+}\times Y$ then from e.g. (\cite{Hale}) and (\ref{Lip}) we have the existence, uniqueness, positivity and continuous of solution to problem (\ref{A2})-(\ref{A22}) in $\mathbb{R}^{+}\times Y.$

\begin{Proposition}\label{Prop1}
There exists  an $M>0$ such that for any solution of (\ref{A2})-(\ref{A22}) there exists a $T>0$ such that,
\begin{equation}\label{M}
S(t)\leq M,\; E(t)\leq M,\; ||E_t||\leq M,\;\ \mbox{and}\; J(t)\leq M   \;\ \mbox{for all}\;\ t\geq T.
\end{equation}

 Moreover,
 \begin{equation}\label{SSS}
  \Lambda \leq  \liminf_{t\rightarrow \infty} S(t),
 \end{equation}
 with  $\Lambda:=\dfrac{A}{\mu+L}.$
\end{Proposition}
\dem
First of all, a simple computation, gives us

\begin{equation*}
\limsup_{t\rightarrow \infty}N(t)\leq \bar{N}.
\end{equation*}
as $S+E+I+R=N$ then
\begin{equation*}
\limsup_{t\rightarrow \infty}(S(t)+E(t))\leq \bar{N}.
\end{equation*}
Therefore there exists $T>0$ such that for all $t\geq T$ we have
\begin{equation*}
S(t)\leq \bar{N},\;\ I(t)\leq \bar{N}\;\ \mbox{and}\;\ E(t)\leq \bar{N}.
\end{equation*}
In addition
\begin{eqnarray*}
||E_t||&=&\sup_{\theta\leq 0} E_t(\theta)=\sup_{u\leq t} E(u)e^{\Delta u}e^{-\Delta t}\\\\
&\leq& \max \{{e^{-\Delta t}||\phi||,Ke^{\Delta T}e^{-\Delta t}, \bar{N}}\},
\end{eqnarray*}
where $K=\sup_{0\leq u\leq T} E(u).$ Further
\begin{eqnarray*}
J(t)&\leq& \alpha\int_0^{\infty}\beta(a)e^{\Delta a}e^{-\mu a}e^{-\Delta a}E_t(-a)da,\\\\
&\leq& \dfrac{\alpha||\beta||}{\mu+\inf(\gamma)-\Delta}||E_t||.
\end{eqnarray*}
Consequently, we can choose $M$ so large such that (\ref{M}) is satisfied. Concerning the estimate (\ref{SSS}). We set $\liminf_{t\rightarrow \infty} S(t)=S_{\infty}$ and $\limsup_{t\rightarrow \infty} J(t)=J^{\infty},$ using the fluctuation method (\cite{Smi-Thie}, \cite{Thieme0}), there exists a sequence $t_k$ such that $S'(t_k)\rightarrow 0$ and $\lim_{t_k\rightarrow \infty} S(t_k)=S_{\infty},$ thus
 \begin{equation*}
 0\geq A-\mu S_{\infty}-f(S_{\infty},J^{\infty}),
 \end{equation*}
 due to (\ref{Lip}), we obtain
 \begin{equation*}
 0\geq A-\mu S_{\infty}-LS_{\infty},
 \end{equation*}
 so,
  \begin{equation*}
 S_{\infty}\geq \dfrac{A}{\mu+L}.
 \end{equation*}
 The proposition is proved.
\cqd

\noindent For the model (\ref{A1}), the number $R_0$ of secondary infections produced by a single infected individual \cite{Diekmann} is defined by

\begin{equation}
R_0=\dfrac{\alpha}{\mu+\alpha}\dfrac{\partial f}{\partial J}(\bar{N},0)\int_0^{\infty}\beta(a)e^{-\mu a}e^{-\int_0^a\gamma(\sigma)d\sigma}da,
\end{equation}

The paper is organized as follows: The next section focused on proving existence of compact attractor of solution of model (\ref{A2})-(\ref{A22}). Section 3 is devoted to prove  that the disease-free equilibrium is globally asymptotically stable whenever $R_0\leq 1.$  Finally, we will investigate the global dynamic of the endemic equilibrium, whenever it exists.

\section{Global compact attractor}

First, it is not difficult to show the existence of a continuous semiflow
\begin{equation}\label{phi0}
 \Phi(t,(S_0,\phi))=(S(t),E_t(.)),
  \end{equation}
  with $(S,E_t)$ is solution of the autonomous problem (\ref{A2})-(\ref{A22}).\\

We choose $X=\mathbb{R}^{+}\times Y$ and we will prove the existence of a compact attractor of all bounded sets of $X$, (the concept of global attractor is presented in e.g. \cite{Magal1},\cite{Smi-Thie}, \cite{Thieme0}).

\begin{Theorem}
 The semiflow $\Phi$ has a compact attractor $\mathbf{A}$ of bounded sets of $X$.
\end{Theorem}
\dem
By Proposition \ref{Prop1}, the semiflow $\Phi$ is point-dissipative.  Hence, by Theorem 2.33 in \cite{Smi-Thie}, we only need to show that  $\Phi$ is eventually bounded on bounded sets and asymptotically smooth in order to prove our Theorem. These two properties are checked by employing the same ideas as in proof of Theorem 6.1 \cite{Rost}.\\

The rest of this section is devoted to describe  some estimates for bounded total trajectories of our system, that are solutions of (\ref{A2})-(\ref{A22}) defined for all $t\in \mathbb{R}.$ These extended solutions play an important role in proving the global asymptotic stability of equilibria see e.g. \cite{Smi-Thie}.\\
\textbf{Total trajectories}\\
 We consider $\bar{\phi}$  a total $\Phi-$trajectory, $\bar{\phi}(t)=(S(t),E_t(.)).$ Then $\bar{\phi}(r+t)=\Phi(t,\bar{\phi}(r)),$ $t\geq 0,$ $r\in \mathbb{R}.$  Thus by a simple computation, a total trajectories satisfy, for all $t\in \mathbb{R},$
 \begin{equation}\label{Total}
 \left\{
 \begin{array}{lll}
 S'(t)=A-\mu S(t)-f(S(t),J(t)),\\\\
 E'(t)=f(S(t),J(t))-(\mu+\alpha)E(t),\\\\
 J(t)=\alpha\int_0^{\infty}\beta(a)e^{-\mu a}\digamma(a) E(t-a)da.
 \end{array}
 \right.
 \end{equation}

 The next lemma provides some useful estimates of the total trajectory when dealing with the compact attractor, $\mathbf{A}.$

\begin{Lemma}\label{lem2}
For all $(S_0,\phi)\in \mathbf{A},$ we have
\begin{equation*}
\begin{array}{lll}
S(t)+E(t)\leq \bar{N},\;\;\;\ \mbox{and}\;\  S(t)\geq \dfrac{A}{\mu+L},\;\ \mbox{where L is the Lipschitz constant,}\\\\
 \mbox{and}\;\ J(t)\leq \dfrac{\alpha\bar{N}}{\mu}||\beta||,
\end{array}
\end{equation*}
for all $t\in \mathbb{R}.$
\end{Lemma}
\dem
First, summing the first and second equations of system (\ref{Total}) we have
 \begin{equation*}
 S'(t)+E'(t)\leq A-\mu(S(t)+E(t)),
 \end{equation*}
 for $t\geq r$ we get
 \begin{equation*}
 (S(t)+E(t))e^{\mu t}\leq (S(r)+E(r))e^{\mu r}+\dfrac{A}{\mu}(e^{\mu t}-e^{\mu r}),
 \end{equation*}
 letting $r\rightarrow -\infty$ we obtain
 \begin{equation*}
 S(t)+E(t)\leq \dfrac{A}{\mu}, \;\ \mbox{for all}\;\ t\in \mathbb{R}.
 \end{equation*}
Also
 \begin{equation*}
 J(t)=\int_0^{\infty}\beta(a)e^{-\mu a}\digamma(a)\alpha E(t-a)da\leq  \dfrac{\alpha A}{\mu^2}||\beta||.
 \end{equation*}

 Now we deal with $S$ in (\ref{Total}). By the boundedness of $J$ and (\ref{Lip}) we have,
\begin{equation*}
\begin{array}{lll}
S'(t)&\geq& A-\mu S(t)-f(S(t),\dfrac{\alpha A}{\mu^2}||\beta||), \\\\
&\geq& A-\mu S(t)-LS(t).
\end{array}
\end{equation*}
Finally, by a straightforward computation,
 \begin{equation}\label{Sestim}
 S(t)\geq \dfrac{A}{\mu+L} \;\;\ \forall  t\in \mathbb{R}.
 \end{equation}
  This complete the proof.
\cqd

\section{The global stability of the disease-free equilibrium}
 This section is devoted to prove the global asymptotic stability of disease-free equilibrium $ (\bar{N}, 0 ).$
 Throughout this section we suppose that the function $f(S,J)$ is concave with respect to $J.$

 \noindent First, remark that system (\ref{A2}) always has a disease-free equilibrium $(\bar{N},0)$.

\begin{Theorem}\label{Free}
 The disease free equilibrium $(\bar{N}, 0 )$ is globally asymptotically stable whenever $R_0\leq 1.$
\end{Theorem}
\dem
Let us introduce the function
\begin{equation*}
\psi(a)=\alpha\dfrac{\partial f}{\partial J}(\bar{N},0)\int_{a}^{\infty}\beta(\sigma)e^{-\mu \sigma}\digamma(\sigma)d\sigma.
\end{equation*}

For $(S_0,\phi)\in \mathbf{A},$ we construct a Lyapunov functional $V(S_0,\phi)=V_1(S_0,\phi)+V_2(S_0,\phi)+\phi(0)$ where

\begin{equation*}
V_1(S_0,\phi)=S_0-\int^{S_0}_{\bar{N}}\lim_{J\rightarrow 0^{+}}\dfrac{f(\bar{N},J)}{f(\eta,J)}d\eta-\bar{N},
\end{equation*}
and
\begin{equation*}
V_2(S_0,\phi)=\int_0^{\infty}\psi(a)\phi(-a)da,
\end{equation*}
Let $\Psi: \mathbb{R}\rightarrow \mathbf{A}$ be a total $\Phi-$trajectory, $\Psi(t)=(S(t),E_t),$ $S(0)=S_0$ and $E_0=\phi,$  with $(S(t),E_t)$ is solution of problem (\ref{Total}).

\begin{equation*}
\begin{array}{lll}
\dfrac{d}{dt}V_1(\Psi(t))&=&(1-\lim_{J\rightarrow 0^{+}}\dfrac{f(\bar{N},J)}{f(S(t),J)})(A-\mu S(t)-f(S(t),J(t)))\\\\
&=& \mu(1-\lim_{J\rightarrow 0^{+}}\dfrac{f(\bar{N},J)}{f(S(t),J)})\big(\bar{N}-S(t)\big)-f(S(t),J(t))\big(1-\lim_{J\rightarrow 0^{+}}\dfrac{f(\bar{N},J)}{f(S(t),J)}\big).
\end{array}
\end{equation*}
 Concerning $V_2$  we have
\begin{equation*}
\begin{array}{lll}
\dfrac{d}{dt}V_2(\Psi(t))&=&\int_0^{\infty}\psi(a)\dfrac{d}{dt}E(t-a)da,\\\\
&=&-\int_0^{\infty}\psi(a)\dfrac{d}{da}E(t-a)da,\\\\
&=&\psi(0)E(t)+\int_0^{\infty}\psi'(a)E(t-a)da\\\\
&=& R_0(\mu+\alpha) E(t)-\alpha\dfrac{\partial f}{\partial J}(\bar{N},0)\int_0^{\infty}\beta(a)e^{-\mu a}\digamma(a)E(t-a)da.
\end{array}
\end{equation*}

Adding $V'_1,$  $V'_2,$ and $E(t)$  gives
\begin{equation*}
\begin{array}{lll}
\dfrac{d}{dt}V(\Psi(t))&=& \mu(1-\lim_{J\rightarrow 0^{+}}\dfrac{f(\bar{N},J)}{f(S(t),J)})\big(\bar{N}-S(t)\big)-\big(1-R_0\big)(\mu+\alpha) E(t)\\\\
&+&f(S(t),J(t))\lim_{J\rightarrow 0^{+}}\dfrac{f(\bar{N},J)}{f(S(t),J)}-\dfrac{\partial f(\bar{N},0)}{\partial J}J(t).
\end{array}
\end{equation*}
 Observe that, the first two terms of this above equation are negative. We claim that, the third term is also negative. Indeed, the concavity of $f$ with respect to $J$ ensures that
\begin{equation*}
f(S,J)\leq J\dfrac{\partial f}{\partial J}(S,0).
\end{equation*}
Hence
\begin{eqnarray*}
f(S(t),J(t))\lim_{J\rightarrow 0^{+}}\dfrac{f(\bar{N},J)}{f(S(t),J)}-J(t)\dfrac{\partial f}{\partial J}(\bar{N},0)&=&f(S(t),J(t))\dfrac{\dfrac{\partial f}{\partial J}(\bar{N},0)}{\dfrac{\partial f}{\partial J}(S(t),0)}-J(t)\dfrac{\partial f}{\partial J}(\bar{N},0),\\\\
&=&\dfrac{\dfrac{\partial f}{\partial J}(\bar{N},0)}{\dfrac{\partial f}{\partial J}(S(t),0)}\big(f(S(t),J(t))-J(t)\dfrac{\partial f}{\partial J}(S(t),0)\big)\leq 0.
\end{eqnarray*}
Notice that,  $\dfrac{d}{dt}V(\Psi(t))=0$ implies that $S(t)=\bar{N}$. Let $Q$ be the largest invariant set, for which  $\dfrac{d}{dt}V(\Psi(t))=0.$ Then in $Q$ we must have $S(t)=\bar{N}$ for all $t\in \mathbb{R}.$ We substitute this into the equation of $S$, we get $J(t)=0$  for all $t\in\mathbb{R}.$ From system (\ref{Total}) we easily get $E(t)=0$ for all $t\in \mathbb{R}.$ Since $\mathbf{A}$ is compact, the $\omega(x)$ and $\alpha(x)$  are non-empty, compact, invariant and attract $\Psi(t)$ as $t\rightarrow \pm \infty,$ respectively. Since $V(\Psi(t))$ is decreasing function of $t,$ $V$ is constant on the $\omega(x)$ and $\alpha(x),$ and thus  $\omega(x)=\alpha(x)=\{(\bar{N},0)\}.$ Consequently $\limit_{t\longrightarrow\pm \infty}\Psi(t)=(\bar{N},0)$ and
\begin{equation*}
\limit_{t\longrightarrow - \infty}V(\Psi(t))=\limit_{t\longrightarrow + \infty}V(\Psi(t))=V(\bar{N},0).
\end{equation*}
 Hence, we obtain $V(\Psi(t))=V(\bar{N},0)$ for all $t\in \mathbb{R}.$ Since $\alpha(x)=\{(\bar{N},0)\}$ then $V(\Psi(t))\leq V(\bar{N},0)$ for all  $t\in \mathbb{R}.$  Since $\Psi$ achieves its minimum value at $(\bar{N},0),$  we conclude that  $\Psi(t)=(\bar{N},0)$ for all $t\in \mathbb{R}.$   In particular $(S_0,\phi)=(\bar{N},0).$ Therefore the attractor $\mathbf{A},$ is the singleton set formed by the disease free equilibrium $(\bar{N},0).$  By Theorem 2.39 in \cite{Smi-Thie}, the disease free equilibrium is globally asymptotically stable.
\cqd
\section{Existence of endemic equilibrium states and Uniform persistence}
In this section, we first ensures the existence of a positive equilibrium states and next, we establish the strongly uniform  persistence of the solution of problem (\ref{Total}).

\begin{Lemma}\label{lem0}

Let $\lim_{J\rightarrow 0^{+}}\dfrac{f(\bar{N},J)}{f(S,J)}>1$ for $S\in [0,\bar{N}).$ Then, if $R_0>1,$  system (\ref{A2}) has  positive endemic equilibrium states.

\end{Lemma}
\dem
An endemic equilibrium is a fixed point of the semiflow $\Phi$,
\begin{equation*}
\Phi(t,(S^{*},E^{*}))=(S^{*},E^{*}), \;\ \mbox{with} \;\ E^{*}\neq 0,\;\ \forall t\geq0.
\end{equation*}

where  $S^{*}, E^{*}$ satisfy
\begin{equation}\label{6}
\left\{
  \begin{array}{lll}
A-\mu S^{*}-f(S^{*},J^{*})=0, \\\\
(\mu+\alpha)E^{*}=f(S^{*},J^{*}),\\\\
J^{*}=\alpha E^{*}\int_0^{\infty}\beta(a)e^{-\mu a}\digamma(a)da.
\end{array}
  \right.
\end{equation}

Combining the equations of  (\ref{6}),  we obtain,
\begin{equation*}
\left\{
  \begin{array}{lll}
A=\mu S^{*}+f(S^{*},J^{*}), \\\\
\bar{D}f(S^{*},J^{*})=J^{*},
\end{array}
  \right.
\end{equation*}
with
\begin{equation}\label{D}
\bar{D}:=\dfrac{\alpha}{\mu+\alpha}\int_0^{\infty}\b(a)e^{-\mu a}\digamma(a)da.
\end{equation}

Following the same arguments as \cite{Feng}, \cite{Korobeinikov}, \cite{Korobeinikov1},  we prove the existence of the positive equilibrium state.
\cqd
We emphasis now on the uniform persistence; For this purpose, we apply Theorem 5.2 in \cite{Smi-Thie},  see also \cite{Hale0}, \cite{Magal1}, \cite{Thieme} for further results in this direction. \\
\noindent We first  make the following assumptions on the incidence $f.$ \\
\noindent We suppose that there exists a positive equilibrium $(S^{*},J^{*})$ verifying (\ref{6}) such that
\begin{equation}\label{S1}
  \left\{
  \begin{array}{lll}
  \dfrac{x}{J^{*}}< \dfrac{f(S,x)}{f(S,J^{*})}< 1 \;\ \mbox{for}\;\ x< J^{*},\\\\
  1< \dfrac{f(S,x)}{f(S,J^{*})}< \dfrac{x}{J^{*}} \;\ \mbox{for}\;\ x> J^{*}.
  \end{array}
  \right.
  \end{equation}
There exists $\e>0$ and there exists $\eta>0$ such that for all $S\in[\bar{N}-\e,\bar{N}+\e]$ we have
\begin{equation}\label{decrea}
  \dfrac{f(S,J_1)}{J_1}\geq \dfrac{f(S,J_2)}{J_2},
\end{equation}
for all $0< J_1 \leq J_2\leq \eta.$

Finally we suppose that

\begin{equation}\label{Pi}
\phi(a)>0,\;\ \mbox{for some} \;\ a\leq 0.
\end{equation}

  \begin{remark}
If the function $f$ is differentiable and concave with respect to $J$ then (\ref{S1}) and (\ref{decrea}) are clearly verified.
\end{remark}

 We define a persistence function $\rho: \mathbb{R}^{+}\times Y_{+}\rightarrow \mathbb{R}^{+}$ by
\begin{equation*}
\rho(S_0,\phi)=E(0),
\end{equation*}
then for $x=(S_0,\phi),$
\begin{equation*}
\rho(\Phi(t,x))=E(t).
\end{equation*}
The following Lemma affirm that the hypothesis $(H1)$ in Theorem 5.2 (\cite{Smi-Thie}) holds.
\begin{Lemma}
Under the hypothesis (\ref{Pi}), the function $\rho$ is positive everywhere on $\mathbb{R}.$
\end{Lemma}
\dem
We first suppose that there exists $r\in \mathbb{R}$ such that $E(t)=0$ for all $t\leq r.$ We claim that $E(t)=0$ for all $t>r.$ Indeed for $t>r$, from the equation of $E$ in (\ref{Total}) and after a change of variable we get
\begin{eqnarray*}
(E(t)e^{(\mu+\alpha)t})'&=&e^{(\mu+\alpha)t} f\big(S(t),\alpha \int_{-\infty}^{t}\beta(t-\sigma)e^{-\mu(t-\sigma)}\digamma(t-\sigma)E(\sigma)d\sigma\big),\\\\
&=& e^{(\mu+\alpha)t}f\big(S(t),\alpha\int_{r}^{t}\beta(t-\sigma)e^{-\mu(t-\sigma)}\digamma(t-\sigma)E(\sigma)d\sigma\big),
\end{eqnarray*}
in view of Lemma \ref{lem2} and (\ref{Lip}) we get
\begin{eqnarray*}
(E(t)e^{(\mu+\alpha)t})'\leq L\bar{N}e^{(\mu+\alpha)t}\alpha\int_{r}^{t}\beta(t-\sigma)e^{-\mu(t-\sigma)}\digamma(t-\sigma)E(\sigma)d\sigma,
\end{eqnarray*}
by integration and Fubini's theorem we have

\begin{eqnarray*}
E(t)e^{(\mu+\alpha)t}&\leq& \alpha L\bar{N} ||\beta||\int_r^{t}e^{(\mu+\alpha)\theta}\int_r^{\theta}e^{-\mu(\theta-\sigma)}E(\sigma)d\sigma d\theta,\\\\
&\leq& \alpha L\bar{N} ||\beta||\int_r^{t}e^{\mu \sigma}E(\sigma)\int_{\theta}^{t}e^{\alpha \theta}d\theta d\sigma,
\end{eqnarray*}
thus, by a straightforward calculation, we obtain
\begin{equation*}
E(t)\leq L \bar{N} ||\beta||\int_r^tE(\sigma)d\sigma.
\end{equation*}

Using Gronwall Lemma and the fact that $E(r)=0$ we conclude that $E(t)=0$ for all $t>r$. This is a contradiction with (\ref{Pi}). Now suppose that there exists a sequence $t_n$ such that $t_n\rightarrow -\infty$ and $E(t_n)>0.$ We set  $E_n(t)=E(t+t_n)$ and $S_n(t)=S(t+t_n).$ So, from Lemma \ref{lem2} and by a simple computation we get

\begin{eqnarray*}
E_n(t)e^{(\mu+\alpha)t}&\geq& E_n(0)+ \int_0^{t}e^{(\mu+\alpha)\theta}f\big( \dfrac{A}{\mu+L}, \alpha\int_{-\infty}^{\theta}\beta(\theta-\sigma)e^{-\mu(\theta-\sigma)}E_n(\sigma)d\sigma \big)d\theta,
\end{eqnarray*}
 since $E_n(0)=E(t_n)>0$ then $E_n(t)>0$ for all $t\geq 0.$ Finally, since $t_n\rightarrow -\infty$ as $n\rightarrow \infty$ then

\begin{eqnarray*}
E(t)>0, \;\ \forall t\in \mathbb{R}.
\end{eqnarray*}
The lemma is proved.
\cqd

As the hypothesis (H1) in Section 5.1. \cite{Smi-Thie}  holds. So from Theorem 5.2. in \cite{Smi-Thie}
it is sufficient to prove the weak uniform persistence in order to establish the strong uniform persistence.

\begin{Theorem}\label{Th1}
Assume that  (\ref{S1}), (\ref{decrea}), (\ref{Pi}) hold. Then there exists some $\eta>0$ such that
\begin{equation*}
\liminf_{t\rightarrow \infty} E(t)>\eta,
\end{equation*}
for all nonnegative solutions of (\ref{A2}) provided that  $R_0>1.$
\end{Theorem}
\dem
By contradiction, we suppose that
\begin{equation*}
\lim_{t\rightarrow \infty}E(t)=0,
\end{equation*}
then we have
\begin{equation*}
\lim_{t\rightarrow \infty}J(t):=\lim_{t\rightarrow \infty}\int_0^{\infty}\alpha \beta(a)e^{-\mu a}\digamma(a)E(t-a)da=0.
\end{equation*}
Now, in view of the equation of $S$ in (\ref{A1}) and for $t$ so large, we have
\begin{eqnarray*}
S'(t)&\geq& A-\mu S(t)-f(\bar{N},J),\\\\
&\geq& A-\mu S(t)-f(\bar{N},\e),
 \end{eqnarray*}
by straightforward computation we obtain
\begin{eqnarray*}
S(t)\geq (S_0-\dfrac{A-f(\bar{N},\e)}{\mu})e^{-\mu t}+\dfrac{A-f(\bar{N},\e)}{\mu},
 \end{eqnarray*}
 and thus
\begin{eqnarray*}
\liminf_{t\rightarrow \infty} S(t)\geq \bar{N}-\psi(\e),
  \end{eqnarray*}
with $\psi(\e)=\dfrac{f(\bar{N},\e)}{\mu}.$
We also have
\begin{equation*}
\limsup_{t\rightarrow \infty}S(t)\leq \bar{N},
\end{equation*}
therefore we can assume that
\begin{equation*}
|S(t)-\bar{N}|<\psi(\e),\;\ t\geq T.
\end{equation*}
Otherwise, since $R_0>1$ then for $\e$ so small and $t^{*}$ so large we have
\begin{equation*}
\dfrac{f(\bar{N}-\psi(\e),\e)}{\e}\int_0^{t^{*}}\alpha\beta(a)e^{-\mu a}\digamma(a)da>\mu+\alpha.
\end{equation*}
Now using the equation of $E$

\begin{equation*}
\begin{array}{lll}
E'(t)&=&f(S(t),J(t))-(\mu+\alpha)E(t),\\\\
&\geq& f(\bar{N}-\psi(\e),J(t))-(\mu+\alpha)E(t),
\end{array}
\end{equation*}
In view of (\ref{decrea})  and the fact that $J(t)<\e$ for $t$ large,  we have,
 $$f(\bar{N}-\psi(\e),J(t))\geq \dfrac{f(\bar{N}-\psi(\e),\e)}{\e}J(t),$$  thus
 \begin{equation*}
\begin{array}{lll}
E'(t)&\geq& \dfrac{f(\bar{N}-\psi(\e),\e)}{\e}J(t)-(\mu+\alpha)E(t),\\\\
&\geq&  \dfrac{f(\bar{N}-\psi(\e),\e)}{\e}\int_0^{t^{*}}\alpha\beta(a)e^{-\mu a}\digamma(a)E(t-a)da-(\mu+\alpha)E(t).
\end{array}
\end{equation*}
Now by employing the same argument as in the proof of Theorem 6.1. (\cite{Rost}) we reach a contradiction.
Hence we conclude.
\cqd

Let $X_0$ be a subset defined as
\begin{equation*}
X_0=\{ (S_0,\phi)\in X ; \phi(a)=0,\;\ \mbox{for all}\;\ a\leq 0  \},
\end{equation*}

\noindent from Theorem 5.7 in \cite{Smi-Thie} we have the following result
\begin{Theorem}\label{Th2}
There exists a compact attractor $\mathbf{A_1}$ that attracts all solutions with initial condition belonging to $X\setminus X_0.$ Moreover $\mathbf{A_1}$ is $\rho-$ uniformly positive, i.e., there exists some $\delta>0$ such that,
\begin{equation}\label{del}
E(t)\geq \delta, \;\  \forall t\in \mathbb{R},\;\ \mbox{and}\;\ \forall  (S_0,\phi)\in \mathbf{A_1}.
\end{equation}
\end{Theorem}
\section{The global stability and uniqueness of the endemic steady state }

In this section, we discuss the global stability of the endemic equilibrium $E^{*}$ of the system (\ref{Total}).
Before stating our main result of this section, we need the following estimate, which guarantees that all solutions of (\ref{Total}) with initial data satisfying (\ref{Pi}), are bounded away from 0.
\begin{Corollary}
There exists $\bar{\delta}>0$ such that, for all $(S_0,\phi)\in \mathbf{A_1},$
\begin{equation*}
J(t):=\alpha\int_0^{\infty}\beta(a)e^{-\mu a}\digamma(a)E(t-a)da \geq \bar{\delta}, \;\
\end{equation*}
and
\begin{equation}\label{f}
f(S(t),J(t))\geq f(\dfrac{A}{\mu+L},\bar{\delta}),
\end{equation}
for all  $t\in\mathbb{R},$ and $\bar{\delta}:=\delta \alpha\int_0^{\infty}\beta(a)e^{-\mu a}\digamma(a)da.$
\end{Corollary}
\dem
Since $\mathbf{A_1}$ is invariant, there exists a total trajectory $\Psi:\mathbb{R}\rightarrow \mathbf{A_1},$ $\Psi(t)=(S(t),E_t)$ with $S(0)=S_0$ and $E_0(a)=\phi(a)$ for $a\leq 0.$ In view of the estimation (\ref{del}),  for all $t\in \mathbb{R}$ we have,
\begin{equation*}
J(t):=\alpha\int_0^{\infty}\beta(a)e^{-\mu a}\digamma(a)E(t-a)da\geq \delta \alpha\int_0^{\infty}\beta(a)e^{-\mu a}\digamma(a)da,
\end{equation*}
Form Lemma \ref{lem2} we easily get (\ref{f}).
\cqd

\begin{Theorem}
 Under the assumptions of Theorem \ref{Th1}, the positive endemic $(S^{*},E^{*})$ is unique and is globally asymptotically stable in $\mathbf{A_1}$.
\end{Theorem}
\dem

Let $\Psi: \mathbb{R}\rightarrow \mathbf{A_1}$ be a total $\Phi-$trajectory, $\Psi(t)=(S(t),E_t(.)),$ $S(0)=S_0$ and $E_0(.)=\phi,$  with $(S(t),E_t(.))$ is solution of problem (\ref{Total}). We define $H(y)=y-ln(y)-1,$ and
\begin{equation}\label{fi}
\psi(a)=f(S^{*},J^{*})\int_a^{\infty}dm(\sigma),\;\  \mbox{with} \;\ dm(\sigma)=\dfrac{\beta(\sigma)e^{-\mu \sigma}\digamma(\sigma)}{\bar{D}}d\sigma.
\end{equation}
with $\bar{D}:=\int_0^{\infty}\beta(\sigma)e^{-\mu \sigma}\digamma(\sigma)d\sigma.$
Notice that $\int_0^{\infty}dm(\sigma)=1.$

\noindent For $(S_0,\phi)\in A_1,$ we consider the following Lyapunov function $V(S_0,\phi)=V_1(S_0,\phi)+V_2(S_0,\phi)+E^{*}V_3(S_0,\phi)$ with
\begin{equation*}
\begin{array}{lll}
V_1(S_0,\phi)=S_0-\int_{S^{*}}^{S_0}\dfrac{f(S^{*},J^{*})}{f(\eta,J^{*})}d\eta-S^{*},
\end{array}
\end{equation*}
and
\begin{equation*}
\begin{array}{lll}
V_2(S_0,\phi)=\int_0^{\infty}\psi(a)H(\dfrac{\phi(-a)}{E^{*}})da.
\end{array}
\end{equation*}
\begin{equation*}
\begin{array}{lll}
V_3(S_0,\phi)=H(\dfrac{\phi(0)}{E^{*}}),
\end{array}
\end{equation*}
Firstly,
\begin{equation}\label{G1}
\begin{array}{lll}
\dfrac{d}{dt}V_1(\Psi(t))=\mu(1-\dfrac{f(S^{*},J^{*})}{f(S(t),J^{*})})(S^{*}-S(t))+(1-\dfrac{f(S^{*},J^{*})}{f(S(t),J^{*})})f(S^{*},J^{*})-(1-\dfrac{f(S^{*},J^{*})}{f(S(t),J^{*})})f(S(t),J(t)).
\end{array}
\end{equation}
Concerning $V_2$ we have
\begin{equation}\label{G2}
\begin{array}{lll}
\dfrac{d}{dt}V_2(\Psi(t))&=&H(\dfrac{E(t)}{E^{*}})\psi(0)+\int_0^{\infty}H(\dfrac{E(t-a)}{E^{*}})\psi'(a)da,\\\\
&=&f(S^{*},J^{*})H(\dfrac{E(t)}{E^{*}})+\int_0^{\infty}H(\dfrac{E(t-a)}{E^{*}})\psi'(a)da,
\end{array}
\end{equation}
Now
\begin{equation}\label{G22}
\begin{array}{lll}
E^{*}V'_3(\Psi(t))=(1-\dfrac{E^{*}}{E(t)})\big(f(S(t),J(t))-(\mu+\alpha)E(t)\big),
\end{array}
\end{equation}
adding $V'_1,$ $V'_2$ and $E^{*}V'_3$ we obtain
\begin{equation*}
\begin{array}{lll}
V'(\Psi(t))&=&\mu(1-\dfrac{f(S^{*},J^{*})}{f(S(t),J^{*})})(S^{*}-S(t))+f(S^{*},J^{*})(1-\dfrac{f(S^{*},J^{*})}{f(S(t),J^{*})})\\\\
&-&f(S(t),J(t))(1-\dfrac{f(S^{*},J^{*})}{f(S(t),J^{*})})\\\\
&+&f(S^{*},J^{*})\big(\dfrac{E(t)}{E^{*}}-\ln\dfrac{E(t)}{E^{*}}-1\big)+\big(1-\dfrac{E^{*}}{E(t)})f(S(t),J(t))\\\\
&-&(\mu+\alpha)E(t)(1-\dfrac{E^{*}}{E(t)})+\int_0^{\infty}H(\dfrac{E(t-a)}{E^{*}})\psi'(a)da.
\end{array}
\end{equation*}
Recall that $(\mu+\alpha)E^{*}=f(S^{*},J^{*}),$ and reorganizing our terms, we have
\begin{equation*}
\begin{array}{lll}
V'(\Psi(t))&=&\mu(1-\dfrac{f(S^{*},J^{*})}{f(S(t),J^{*})})(S^{*}-S(t))+f(S(t),J(t))\dfrac{f(S^{*},J^{*})}{f(S(t),J^{*})}+\int_0^{\infty}H(\dfrac{E(t-a)}{E^{*}})\psi'(a)da\\\\
&+&f(S^{*},J^{*})(1-\dfrac{f(S^{*},J^{*})}{f(S(t),J^{*})})-f(S^{*},J^{*})\ln\dfrac{E(t)}{E^{*}}-\dfrac{E^{*}}{E(t)}f(S(t),J(t)).
\end{array}
\end{equation*}
Adding and subtracting the term
\begin{equation*}
f(S^{*},J^{*})\big(\ln\dfrac{f(S(t),J(t))}{f(S(t),J^{*})}+1\big),
\end{equation*}
and using the definition of the function $H$ we get

\begin{equation*}
\begin{array}{lll}
V'(\Psi(t))&=&\mu(1-\dfrac{f(S^{*},J^{*})}{f(S(t),J^{*})})(S^{*}-S(t))+\int_0^{\infty}H(\dfrac{E(t-a)}{E^{*}})\psi'(a)da+f(S^{*},J^{*})H(\dfrac{f(S(t),J(t))}{f(S(t),J^{*})})\\\\
&+&f(S^{*},J^{*})\big(\ln\dfrac{f(S(t),J(t))}{f(S(t),J^{*})}-\ln\dfrac{E(t)}{E^{*}}-\dfrac{E^{*}f(S(t),J(t))}{E(t)f(S^{*},J^{*})}-\dfrac{f(S^{*},J^{*})}{f(S(t),J^{*})}+2\big),
\end{array}
\end{equation*}
from this and the fact that,
\begin{equation*}
\ln\dfrac{f(S(t),J(t))}{f(S(t),J^{*})}=\ln\dfrac{f(S(t),J(t))}{f(S^{*},J^{*})}+\ln\dfrac{f(S^{*},J^{*})}{f(S(t),J^{*})},
\end{equation*}
 we find,
\begin{equation*}
\begin{array}{lll}
V'(\Psi(t))&=&\mu(1-\dfrac{f(S^{*},J^{*})}{f(S(t),J^{*})})(S^{*}-S(t))+\int_0^{\infty}H(\dfrac{E(t-a)}{E^{*}})\psi'(a)da+f(S^{*},J^{*})H(\dfrac{f(S(t),J(t))}{f(S(t),J^{*})})\\\\
&+&f(S^{*},J^{*})\big(\ln\dfrac{E^{*}f(S(t),J(t))}{E(t)f(S^{*},J^{*})}-\dfrac{E^{*}f(S(t),J(t))}{E(t)f(S^{*},J^{*})}+\ln\dfrac{f(S^{*},J^{*})}{f(S(t),J^{*})}-\dfrac{f(S^{*},J^{*})}{f(S(t),J^{*})}+2\big),
\end{array}
\end{equation*}

therefore, in view of the definition of the function $\psi,$
\begin{equation*}
\begin{array}{lll}
V'(\Psi(t))&=&\mu(1-\dfrac{f(S^{*},J^{*})}{f(S(t),J^{*})})(S^{*}-S(t))+f(S^{*},J^{*})\big(H(\dfrac{f(S(t),J(t))}{f(S(t),J^{*})})-\int_0^{\infty}H(\dfrac{E(t-a)}{E^{*}})dm(a)da\big)\\\\
&-&f(S^{*},J^{*})\big\{H(\dfrac{E^{*}f(S(t),J(t))}{E(t)f(S^{*},J^{*})})+H(\dfrac{f(S^{*},J^{*})}{f(S(t),J^{*})})\big\}.
\end{array}
\end{equation*}

Since the function $f$ is nondecreasing with respect to $S$, the first term is negative. Now we set
\begin{equation*}
X=H(\dfrac{f(S(t),J(t))}{f(S(t),J^{*})})-\int_0^{\infty}H(\dfrac{E(t-a)}{E^{*}})dm(a),
\end{equation*}
and we claim that $X$ is negative. Indeed,  using the fact that $H$ is convex, and by Jensen inequality see e.g. \cite{Lieb}, \cite{Thieme0} we get
\begin{eqnarray*}
X &\leq& H(\dfrac{f(S(t),J(t))}{f(S(t),J^{*})})-H\big(\int_0^{\infty}\dfrac{E(t-a)}{E^{*}}dm(a)\big),
\end{eqnarray*}
from the definition of $dm(a)$ and $J$ we have,
\begin{eqnarray*}
X &\leq& H(\dfrac{f(S(t),J(t))}{f(S(t),J^{*})})-H\big(\dfrac{J(t)}{J^{*}}\big).
\end{eqnarray*}
\noindent Let us consider a time $t$ where $Z:=\dfrac{J(t)}{J^{*}}<1;$\\

according to the hypothesis (\ref{S1})  we obtain
\begin{eqnarray*}
\dfrac{f(S(t),J(t))}{f(S(t),J^{*})}\geq \dfrac{J(t)}{J^{*}}.
\end{eqnarray*}

Hence, from $H(1)=0,$  $H$ is decreasing in $(0,1),$ and $f$ is nondecreasing with respect to $J$ we get,

\begin{equation*}
H\big(\dfrac{f(S(t),J(t))}{f(S(t),J^{*})}\big)\leq H(\dfrac{J(t)}{J^{*}}),
\end{equation*}
and thus $X\leq 0.$

For  other values of $t,$  i.e., $Z>1,$ again from (\ref{S1}), we have
\begin{eqnarray*}
\dfrac{f(S(t),J(t))}{f(S(t),J^{*})}\leq \dfrac{J(t)}{J^{*}}.
\end{eqnarray*}

and thus, ($H$ is increasing in $(1,\infty)$)
\begin{equation*}
H(\dfrac{f(S(t),J(t))}{f(S(t),J^{*})})\leq H(\dfrac{J(t)}{J^{*}}).
\end{equation*}

Consequently the claim is proved and so $\dfrac{dV}{dt}\leq 0.$ \\

\noindent Notice that $\dfrac{d}{dt}V(\Psi(t))=0$ implies that $S(t)=S^{*}.$ Let $Q$ be the largest invariant set, for which  $\dfrac{d}{dt}V(\Psi(t))=0.$ Then in $Q$ we must have $S(t)=S^{*}$ for all $t\in \mathbb{R}$ and

\begin{equation*}
H(\dfrac{f(S(t),J(t))E^{*}}{f(S(t),J^{*})E(t)})=0.
\end{equation*}
and thus
\begin{equation}\label{eq}
f(S^{*},J(t))E^{*}=f(S^{*},J^{*})E(t)
\end{equation}
Using the equation of $S$ in (\ref{Total}) we obtain
\begin{equation*}
A-\mu S^{*}=f(S^{*},J(t)),
\end{equation*}

from this, and (\ref{6}) we have
\begin{equation*}
f(S^{*},J(t))=f(S^{*},J^{*}),
\end{equation*}
 then $J(t)=J^{*}$ for all $t\in \mathbb{R}.$ Substituting this result in (\ref{eq}), we get
\begin{equation*}
E(t)=E^{*}, \;\ \mbox{for all} \;\ t\in \mathbb{R}.
\end{equation*}
 Following the same arguments as in proof of Theorem \ref{Free} we conclude the global asymptotic stability of the endemic equilibrium.  As the equality $\frac{d}{dt}V(\Psi(t))=0$ only holds on the line $S=S^{*},$ the endemic equilibrium is unique.
\cqd

\end{document}